\documentclass[11pt,a4paper]{article}
\usepackage[numbers, sort & compress]{natbib}
\usepackage{lmodern}
\usepackage{amssymb}
\usepackage{pstricks}
\usepackage{pst-node}
\usepackage{lineno}
\usepackage[normalem]{ulem}
\usepackage{amsmath}
\usepackage{amsthm}
\usepackage{amsfonts}
\usepackage{tcolorbox}
\usepackage{multicol}

\newtheorem{thm}{Theorem}[section]
\newtheorem{lem}[thm]{Lemma}

\newtheorem{cor}[thm]{Corollary}

\renewcommand\caption[1]{\small\refstepcounter{figure}%
\begin{center}\textbf{Fig.\ \thefigure .}\ #1\end{center}\normalsize}

\newrgbcolor{gray5}{0.5 0.5 0.5}
\newrgbcolor{gray7}{0.7 0.7 0.7}
\newrgbcolor{gray9}{0.9 0.9 0.9}
\newrgbcolor{gray95}{0.95 0.95 0.95}
\newrgbcolor{gray85}{0.85 0.85 0.85}
\newrgbcolor{gray75}{0.75 0.75 0.75}
\newrgbcolor{gray8}{0.8 0.8 0.8}

\def\G{{\ensuremath{G}}}
\def\T{{\ensuremath{T}}}

\def\Dpr{\ensuremath{D_{\rm pr}}}
\def\Dev{\ensuremath{D_{\rm ev}}}
\def\gpr{\ensuremath{\gamma_{\rm pr}}}
\def\gev{\ensuremath{\gamma_{\rm ev}}}

\def\Xev{\ensuremath{X_{\rm ev}}}

\def\Yev{\ensuremath{Y_{\rm ev}}}

\def\VG{\ensuremath{V_{\G}}}
\def\EG{\ensuremath{E_{\G}}}
\def\NG{\ensuremath{N_{\G}}}

\def\graphG{\ensuremath{\G=(\VG,\EG)}}

\newrgbcolor{darkblue}{0.2 0.3 0.6}
\newrgbcolor{lightblue}{0.5 0.6 0.9}
\newrgbcolor{darkgreen}{0.4 0.7 0.4}
\newrgbcolor{darkdarkgreen}{0.2 0.5 0.2}
\newrgbcolor{darkred}{0.7 0.4 0.1}

\usepackage{color}                          % colors
\usepackage{framed}                         % frames
\usepackage{mdframed}                       % better frames

\usepackage{fancyhdr}                   % Fancy headings

\definecolor{prjBgColor}{RGB}{248,248,237}
\definecolor{prjLnColor}{RGB}{180,200,224}

\begin{document}

\begin{center}
\noindent{\LARGE{\textsc{Unique Paired vs Edge-\!Vertex\\[2mm] Minimum Dominating Sets in Trees}}}
\\[8mm]
\small
\textsc{Mateusz Miotk}, 
\textsc{Michał Zakrzewski},
\textsc{Paweł Żyliński}
\\[2mm]
University of Gda\'{n}sk, 80-308 Gda\'nsk, Poland\\[2mm]
\texttt{\{mateusz.miotk, michal.zakrzewski,pawel.zylinski\}@ug.edu.pl}
\end{center}

\vspace{3mm}
\normalsize
\begin{abstract}\small
\noindent{}We prove that the class of trees with unique minimum edge-vertex dominating sets is equivalent to the class of trees with unique minimum paired dominating sets. 
\\[2mm]
\end{abstract}

\section{Introduction}
A subset $D$ of \VG\ is said to be a {\em dominating set} of a graph $\G=(\VG,\EG)$ if each vertex belonging to the set $\VG \setminus D$ has a~neighbour in $D$~\cite{Berge58,Ore62}. A set $D\subseteq \VG$ is a {\em  paired-dominating set\/} of \G\ if $D$ is a dominating set of \G\ and the induced subgraph $\G[D]$ has a perfect matching. The {\em paired-domination number\/} of~\G, denoted by $\gamma_{\rm pr}(\G)$, is defined to be the minimum cardinality of a paired-dominating set $D$ of \G, and any minimum paired dominating set of \G\ is referred to as a $\gamma_{\rm pr}$-set~\cite{HS95,HS98}. %A constructive characterization of trees with unique minimum paired dominating sets can be found in~\cite{CH04}. 
Next, an edge $e \in \EG$ is said to {\em $ev$-dominate\/} a~vertex $v \in \VG$ if $e$ is incident to $v$ or $e$ is incident to a vertex adjacent to $v$. A set $M \subseteq \EG$ is an {\em edge-vertex dominating  set\/} (or simply, an {\em $ev$-dominating set\/}) of~\G\ if each vertex of \G\ is $ev$-dominated by some edge in~$M$. Finally, the {\em edge-vertex domination number of \G\/}, denoted by $\gamma_{ev}(\G)$, is the minimum cardinality of an $ev$-dominating set of~\G, and any minimum edge-vertex dominating set of \G\ is referred to as a $\gamma_{\rm ev}$-set~\cite{LP85}.

Clearly, any matching in a minimum paired-dominating set of a graph $G$ constitutes an edge-vertex dominating set of $G$, but not vice versa, that is, all the vertices of edges forming a minimum edge-vertex dominating set of $G$ do not always constitute a paired-dominating set of~$G$. However, a fundamental, but not common, relation between edge-vertex domination and paired domination was proved by Hedetniemi et al.~\cite{HHH11} who established that for any graph \G\ with no isolated vertex we have $2\gamma_{\rm ev}(\G)=\gamma_{\rm pr}(\G)$. 

\begin{thm}\label{thm:HHH}{\em\cite{HHH11}}   
If \G\ is a graph with no isolated vertex, then $2\gamma_{\rm ev}(\G)=\gamma_{\rm pr}(\G)$.  
\end{thm}

\noindent In this short note, we present another one, when restricted to unique dominating sets. Namely, a commonly used approach for constructive characterizations of trees with unique minimum dominating sets, for different models of domination~\cite{BCLM11,BBZIE20,CH04,CH10,CR14,FV02,GHMR94,HH02,HH20a,
%HH20b,
HH21,H17,L06,LHXL11,Sharada15,SChNV23,SChNY2025,T93,ZWZ18}
is to define a set of basic operations such that starting from the initial pre-defined (finite) set of small trees, playing a role of seeds, any tree with the unique minimum dominating set (in the relevant model) can be constructed by applying, successively, a finite sequence of these operations. 
In particular, Chellali and Haynes~\cite{CH04} characterized the class of trees with unique minimum paired dominating sets by providing the relevant set of four basic graph operations, whereas Senthilkumar et al.~\cite{SChNY2025} -- the relevant set of five basic graph operations in the case of edge-vertex domination. However, we establish the following result.

\begin{thm}\label{thm:main} The class of trees having unique minimum edge-vertex dominating sets is equivalent to the class of trees having unique minimum paired dominating sets. 
\end{thm}

Observe that if a graph \G\ has a cycle, then uniqueness of a \gpr-set of \G\ does not imply uniqueness of a \gev-set of \G, see Figure~\ref{fig:main} for an example. However, as a side result of the proof of Theorem~\ref{thm:main}, we obtain that such a $\gpr$-uniqueness forces all \gev-sets of \G\ to be spanned on the same subset of vertices (Corollary~\ref{cor:general} and~\ref{cor:general2}).     

\begin{figure}
    \centering
    
\begin{center}
\pspicture(0,0)(3,3)
%\psgrid

\cnode[linewidth=0.5pt,linecolor=black,fillstyle=solid,fillcolor=lightgray](1,1){3pt}{c1}
\cnode[linewidth=0.5pt,linecolor=black,fillstyle=solid,fillcolor=lightgray](1,2){3pt}{c2}
\cnode[linewidth=0.5pt,linecolor=black,fillstyle=solid,fillcolor=lightgray](2,2){3pt}{c3}
\cnode[linewidth=0.5pt,linecolor=black,fillstyle=solid,fillcolor=lightgray](2,1){3pt}{c4}

\cnode[linewidth=0.5pt,linecolor=black,fillstyle=solid,fillcolor=white](0.2,0.2){3pt}{x1}
\cnode[linewidth=0.5pt,linecolor=black,fillstyle=solid,fillcolor=white](0.2,2.8){3pt}{x2}
\cnode[linewidth=0.5pt,linecolor=black,fillstyle=solid,fillcolor=white](2.8,2.8){3pt}{x3}
\cnode[linewidth=0.5pt,linecolor=black,fillstyle=solid,fillcolor=white](2.8,0.2){3pt}{x4}

\ncline[linewidth=0.5pt]{-}{x1}{c1}
\ncline[linewidth=0.5pt]{-}{x2}{c2}
\ncline[linewidth=0.5pt]{-}{x3}{c3}
\ncline[linewidth=0.5pt]{-}{x4}{c4}

\ncline[linewidth=0.5pt]{-}{c2}{c1}
\ncline[linewidth=0.5pt]{-}{c2}{c3}
\ncline[linewidth=0.5pt]{-}{c4}{c3}
\ncline[linewidth=0.5pt]{-}{c4}{c1}

\endpspicture
\end{center}
    \caption{The graph \G\ has a unique \gpr-set (gray vertices), whereas its \gev-set is not unique: there are two such sets, but both spanned on the same set of (gray) vertices.}
    \label{fig:main}
\end{figure}

\section{Main Results}\label{sec-2}

We assume that all graphs are simple, i.e., they contain neither loops nor parallel edges. Recall that the {\em open neighborhood} of a vertex $v$ in a connected graph $\G = (\VG,\EG)$ is the set $\NG(v) = \{u \in \VG \colon uv \in \EG\}$, whereas for an edge subset $M \subseteq \EG$, $\VG(M)$ denotes the set of vertices of all edges in $M$. We start with the following crucial lemma (its proof, in general, follows the argument in the proof of Theorem 73 in~\cite{L07} and Theorem 2.2 in~\cite{HHH11}, but proceeds with a more detailed analysis).

\begin{lem}\label{lem:main}
Let \Dev\ be a \gev-set of a graph \graphG. If there exist two 
%$($distinct$)$
edges in $\Dev$ sharing a vertex, then there exist two $($distinct$)$ \gev-sets $\Dev'$ and $\Dev''$ of \G\ such that $\VG(\Dev') \neq \VG(\Dev'')$ and no two edges in $\Dev'$, resp. $\Dev''$, share a vertex. 
\end{lem}

\begin{proof}
Let \Dev\ be a \gev-set of a graph \graphG. We have the following claim (it follows immediately from minimality of \Dev).
\\[2mm]
\textbf{Claim.} {\sl There are no three $($distinct$)$ edges in \Dev\ such that they constitute either a~$4$-vertex path or a $3$-vertex cycle in \G.}
\\[3mm]
Assume now that there exist two edges $e_1,e_2 \in \Dev$ sharing a vertex, say $e_1=x_1x_2$ and $e_2=x_2x_3$. Since \Dev\ is a \gev-set, there exists $x_0 \in \NG(x_1)$ such that $x_0$ is ev-dominated only by edge $e_1$, and analogously, there exists $x_4 \in \NG(x_3)$ such that $x_4$ is ev-dominated only by edge $e_2$. Consider now the sets $\Xev^1=(\Dev\setminus\{e_1\}) \cup\{x_0x_1\}$ and $\Yev^1=(\Dev\setminus\{e_2\}) \cup\{x_3x_4\}$; we shall refer to such edge replacement operation on \Dev\ as {\em \Dev-twinning}. Clearly, $\Xev^1$ and $\Yev^1$ are distinct \gev-sets of \G\ and $\Xev^1 \cap \Yev^1 = \Dev \setminus \{e_1,e_2\}$. Moreover, taking into account the above claim, each of them has less (but the same) number of pairs of edges sharing a vertex. Consequently, if they have no such pair of distinct edges having a vertex in common, we are done. 

Otherwise, we focus only on the set $\Xev^1$ and apply $\Xev^1$-twinning, which now results in two distinct \gev-sets $\Xev^2$ and $\Yev^2$. We continue (iteratively) $\Xev^i$-twinning procedure as long as $\Xev^i$ (and so $\Yev^i$) has at least one pair of edges sharing a vertex. Eventually, we arrive to the situation when $\Xev^i$ (and so $\Yev^i$) has no two edges having a vertex in common. Clearly, keeping in mind minimality of the initial set \Dev, it follows from the construction that the sets $\Xev^i$ and $\Yev^i$ are the two sought $\gev$-sets of \G. \end{proof}

%\noindent As a direct consequence of Lemma~\ref{lem:main}, we have the following corollary.
%, which provides a method to construct such a set.

%\begin{cor}{\em\cite{HHH11}}
%For any graph \G, there always exists a \gev-set of \G\ with no two edges sharing a vertex.  
%\end{cor}

Next, let \G\ be a graph of order at least three (the base case of a $2$-vertex tree is obvious) and assume that all \gev-sets of \G\ are spanned on the same (unique) vertex set~$D$. Let \Dev\ be any \gev-set of \G.
Taking into account Lemma~\ref{lem:main},  uniqueness of $D$ implies that no two elements of \Dev\ share a vertex. Consequently, the set $D$ of size $2|\Dev|$ -- perfectly matchable by \Dev\  -- is a minimum paired-dominating set of \G\ by~Theorem~\ref{thm:HHH}. Furthermore, we claim that $D$ is unique one. Indeed, suppose to the contrary that there exists another \gpr-set\ $D'$ of \G. Then, a perfect matching $M'$ in the induced subgraph $\G[D']$ is an edge-vertex dominating set of \G\, with $|D'|/2=|\Dev|=\gev(\G)$,  which contradicts uniqueness of $D$ (since $D' \neq D$). %Consequently, \Dpr\ is a unique minimum paired-dominating set of \T. 

\begin{cor}
    If a graph \G\ has a unique \gev-set, then \G\ has a unique \gpr-set. 
\end{cor}

Assume now that a graph \G\ has a unique \gpr-set \Dpr. We first claim that $V_\G(\Dev')=V_\G(\Dev'')$ for any two \gev-sets of \G. Indeed, consider a \gev-set \Dev\ of~\G, with $2|\Dev|=|\Dpr|$ (by Theorem~\ref{thm:HHH}), imposed by a perfect matching in $\G[\Dpr]$. It follows from Lemma~\ref{lem:main} that no two edges in (any) \gev-set \Dev\ of~\G\ share a vertex (since otherwise, the relevant sets $V_\G(\Dev')$ and $V_\G(\Dev'')$ in Lemma~\ref{lem:main} are two distinct \gpr-sets of \G, a contradiction with \Dpr\ being unique). So suppose now that \G\ has two distinct minimum edge-vertex dominating sets, say $\Xev'$ and $\Xev''$, with no two edges sharing a vertex. If $V_\G(\Xev') \neq V_\G(\Xev'')$, then both $V_\G(\Xev')$ and $V_\G(\Xev'')$ are paired-dominating sets of \G\ of size $\gpr(\G)$ -- a contradiction with uniqueness of $D_{\rm pr}$. Therefore, we must have $V_\G(\Xev') = V_\G(\Xev'')$, which makes us in a position to prove Theorem~\ref{thm:main}.

Observe that if $V_\G(\Xev') = V_\G(\Xev'')$, then the induced subgraph $\G[V_\G(\Xev') \cup V_\G(\Xev'')]$ has a cycle (since both $\Dev'$ and $\Dev''$ are perfect matchings in $G[V_\G(\Xev') \cup V_\G(\Dev'')]$), and thus \G\ is not a tree. Therefore, if \G\ is a tree, then we must have $\Xev'=\Xev''$, which eventually completes the (simple) proof of Theorem~\ref{thm:main}.\qed 
\\[2mm]
We note in passing that we have actually proved the following two properties.
\begin{cor}\label{cor:general}~
 Let \T\ be a non-trivial tree. If \Dpr\ is a unique \gpr-set of \T, then \T\ has a unique \gev-set which is the $($unique$)$ perfect matching of $V_\T(\Dev)=\Dpr$. Analogously, if~\Dev\ is a unique \gev-set of \T, then  \T\ has a unique \gpr-set which is $V_\T(\Dev)$.
\end{cor}

\begin{cor}\label{cor:general2}~
Let \G\ be a simple graph without isolated vertices. Then \G\ has a unique minimum paired dominating set, say $D$, if and only if all its minimum edge-vertex dominating sets are spanned on all and only vertices of $D$.
\end{cor}

\end{document}